\documentclass{ifacconf}

\usepackage{natbib} 

\usepackage{graphicx} 
\usepackage[caption=false,font=normalsize,labelfont=sf,textfont=sf]{subfig}

\usepackage{amsmath}
\usepackage{amssymb}

\usepackage{nicefrac}

\usepackage[width=\textwidth]{caption}

\usepackage{mathtools}
\usepackage[rightComments=false, indLines=false, beginLComment=/*~, endLComment=*/~]{algpseudocodex}
\usepackage{algorithm}

\graphicspath{{images/}}
\DeclareGraphicsExtensions{.pdf,.jpeg,.png}
\usepackage{xcolor}

\begin{document}
\begin{frontmatter}

    \title{Neural ODEs as Feedback Policies for Nonlinear Optimal Control\thanksref{footnoteinfo}}

    \thanks[footnoteinfo]{This work has been submitted to IFAC for possible publication.}

    \author[First]{Ilya Orson Sandoval}
    \author[First]{Panagiotis Petsagkourakis}
    \author[First]{Ehecatl Antonio del Rio-Chanona}

    \address[First]{Centre for Process Systems Engineering\\ Imperial College London, London, United Kingdom\\ (os220@ic.ac.uk, p.petsagkourakis@imperial.ac.uk \& a.del-rio-chanona@imperial.ac.uk).}

    \begin{abstract}  
        Neural ordinary differential equations (Neural ODEs) define continuous time dynamical systems with neural networks.
        The interest in their application for modelling has sparked recently, spanning hybrid system identification problems and time series analysis.
        In this work we propose the use of a neural control policy capable of satisfying state and control constraints to solve nonlinear optimal control problems.
        The control policy optimization is posed as a Neural ODE problem to efficiently exploit the availability of a dynamical system model.
        We showcase the efficacy of this type of deterministic neural policies in two constrained systems: the controlled Van der Pol system and a bioreactor control problem.
        This approach represents a practical approximation to the intractable closed-loop solution of nonlinear control problems.
    \end{abstract}


\end{frontmatter}

\section{Introduction}

Neural policies represent the dominant approach to parametrize controllers in Reinforcement Learning (RL) research.  
Their attractiveness relies on their dimensional scaling property and universal approximation capacity.  
However, their training procedure usually relies on inefficient sampling based strategies to estimate the gradient of the objective function to be optimized.  
On the contrary, when an environment model is available as a differential equation system, it is possible to leverage it for efficiency through methods based on optimal control and dynamic optimization \citep{ainsworth_FasterPolicy_2021, yildiz_ContinuoustimeModelbased_2021}.
In this work we explore this setting, exploiting a neural policy to parametrize a deterministic control function as a state feedback controller.
This approach provides an approximation to the practically intractable optimal closed-loop policy in continuous time.

The optimization of such a controller follows the same strategy as the training of Neural Ordinary Differential Equations (Neural ODEs) \citep{chen_NeuralOrdinary_2018}.
In our application, the weights of the network only define the control function within a predefined system, instead of defining the whole differential equation system as in NeuralODEs.
To understand the training procedure, it is fruitful to overview the close connection between adjoint sensitivity analysis used in dynamical systems and the backpropagation algorithm used in neural networks.
We revise the literature surrounding these topics and their use in recent applications where both fields meet.

Backpropagation may be seen as a judicious application of the chain rule in computational routines, introduced in optimal control \citep{griewank_reverse_2012} and popularized within the neural network community \citep{{rumelhart1986learning}}.
The first derivations trace back to the introduction of the \textit{Kelley-Bryson gradient method} to solve multistage nonlinear control problems \citep{dreyfus_numerical_1962}.  
In neural networks, it may be derived from the optimization problem where Lagrange multipliers (adjoint variables) enforcing the transitions between states \citep{mizutani_DerivationMLP_2000}.  



On continuous time optimal control problems, the optimality conditions from Pontryagin's Maximum Principle \citep{pontryagin_MathematicalTheory_1986} establish a connection between the sensitivities of a functional cost and the adjoint variables.
This relationship is exploited in continuous sensitivity analysis, where it is used to estimate the influence of parameters in the solution of differential equation systems \citep{serban_cvodes_2005, jorgensen_AdjointSensitivity_2007}.  
When the time is discretized in an ODE system, backpropagation is analogous to the adjoint system of the maximum principle \citep{griewank_reverse_2012, baydin_automatic_2017}, and its use to propagate sensitivities is called discrete sensitivity analysis.  
Modern differential equation solvers include implementations of either continuous or discrete sensitivity analysis, relying on the solution of secondary differential equation systems (\textit{optimize-then-discretize}) or automatic differentiation of the integrator routines (\textit{discretize-then-optimize}). \footnote{
    The \textit{discretize-then-optimize} distinction has a different meaning in optimal control literature \citep{biegler_NonlinearProgramming_2010}.  
    There it refers to \textit{direct approaches} to the optimization problem, which first discretize the dynamical equations to afterwards solve the finite-dimensional nonlinear optimization.
    The \textit{optimize-then-discretize} refers to any procedure that departs from the optimality conditions of the problem, also called \textit{indirect approaches}.
    Since the adjoint equations are part of the optimality conditions in dynamic optimization, the classic \textit{indirect} classification includes both variants in ML literature.
    From an optimization perspective, there is not difference in the approach that calculates the gradients since both do so through adjoints; the selection of the approach is merely practical, based on the peculiarities of each problem \citep{ma_ComparisonAutomatic_2021}.
    \label{fnote:terminology}
}

The adjoint approach suggested with the introduction of Neural ODEs \citep{chen_NeuralOrdinary_2018} is a modern variant of Control Vector Iteration (CVI) \citep{luus_ControlVector_2009}, a sequential indirect strategy that optimizes a fixed vector of parameters by relaxing just one of the necessary conditions of optimality \footnote{This is covered in detail in section \ref{sec:adjoint_sensitivity}.}.  
The relaxed condition is approximated iteratively along optimization rounds where the dynamical and adjoint equations are always satisfied, as in feasible paths methods \citep{chachuat_NonlinearDynamic_2007}.
An algorithmic improvement introduced by \citep{chen_NeuralOrdinary_2018} is the efficient use of reverse-mode Automatic Differentiation (AD) to calculate vector-Jacobians products that appear through the differential equations defining the optimality conditions of the problem.
This grants the possibility of dealing with high dimensional parameter problems efficiently, which is crucial for neural policies within dynamical systems in continuous time.
Furthermore, it crucially avoids the symbolic derivations historically associated with indirect methods including CVI in the numerical optimal control literature \citep{biegler_NonlinearProgramming_2010}.  

The use of neural networks within control systems was explored originally in the 90s \citep{chen_BackpropagationNeural_1989, miller_NeuralNetworks_1990} where the focus was on discrete time systems \citep{hunt_NeuralNetworks_1992}.  
The use of neural control policies in this discrete vein has also attracted attention recently in nonlinear control \citep{rackauckas_GeneralizedPhysicsinformed_2020, adhau_ConstrainedNeural_2021,jin_PontryaginDifferentiable_2020} and MPC \citep{amos_DifferentiableMPC_2018, karg_EfficientRepresentation_2020, drgona_LearningConstrained_2022}.
The gradients required for optimization are computed through direct application of AD over the evolution of a discrete system; a strategy coined as \textit{differentiable control} or more generally \textit{differentiable simulations}.
These strategies revive some of the original ideas that brought interest to neural networks in nonlinear control \citep{cao_FormulationNonlinear_2005} with modern computational tooling for AD calculations \citep{baydin_automatic_2017}.  


Reinforcement Learning (RL) commonly leverage neural parametrization of policies within Markov decision processes (MPD) \citep{sutton_ReinforcementLearning_1992} and has had a rising success with dimensionality scaling thanks to deep learning \citep{schmidhuber_DeepLearning_2015}.
The adaptation of RL approaches to continuous time scenarios based on Hamilton-Jacobi-Bellman formulations was originally explored in \citep{munos_ConvergentReinforcement_1997, doya_ReinforcementLearning_2000}.
A continuous time actor-critic variation based on discrete time data was analysed more recently in \citep{yildiz_ContinuoustimeModelbased_2021}.
The study of the model-free policy gradient method to continuous time was explored by \citep{munos_PolicyGradient_2006}.
In continuous time settings with deterministic dynamics as an environment, neural policies may be trained with techniques from dynamic optimization like CVI, avoiding noisy sampling estimations as in classic RL.
A comparison of the training performance improvement between deterministic model-based neural policies and model-free policy gradient was showcased in \citep{ainsworth_FasterPolicy_2021}.

With a view in practical applications, it is crucial to be able to satisfy constraints while also optimizing the policy performance.
While this has been a major focus in optimal control since its inception \citep{bryson_applied_1975}, there has been little attention to general nonlinear scenarios.
Most works assume either linear dynamics or fixed control profiles instead of state feedback policies.
In relation to model-free RL, constraint enforcement is an active area of research \citep{brunke_SafeLearning_2021}.
Recent work has explored variants of objective penalties \citep{achiam_ConstrainedPolicy_2017a}, Lyapunov functions \citep{chow_LyapunovbasedSafe_2019} and satisfaction in chance techniques \citep{petsagkourakis_ChanceConstrained_2022}.

Here we develop a strategy that allows continuous time policies to solve general nonlinear control problems while satisfying constraints successfully.
Our approach is based on the deterministic calculation of the cost functional gradient with respect to the static feedback policy parameters given a white-box dynamical system environment.
Saturation from the policy architecture enforces hard control constraints while state constraints are enforced through relaxed logarithmic penalties and an adaptive barrier update strategy.
We furthermore showcase how the inclusion of the feedback controller within the ODE definition shapes the whole phase space of the system.
This Neural ODE quality is impossible to achieve in nonlinear systems with standard optimal control methods that only provide controls as a function of time.

\section{Problem statement}

We consider continuous time optimal control problems with fixed initial condition and fixed final time.
The cost functional is in the Bolza form, including both a running ($\ell$) and a terminal cost ($\phi$) in the functional objective ($J$) \citep{bryson_applied_1975}:  
\begin{equation} \label{eq:problem_statement}
    \begin{aligned}
        \min_{\theta} \quad & J = \int_{t_0}^{t_f}{\ell\left(x(t), \pi_\theta(x)\right)} \,dt + \phi(x(t_f)), \\
        \textrm{s.t.} \quad & \dot x(t) = f(x(t), \pi_\theta(x), t),                                          \\
                            & x(t_0) = x_0,                                                                   \\
                            & g(x(t), \pi_\theta(x)) \leq 0,                                                  \\
    \end{aligned}
\end{equation}
where the time window is fixed $t \in [t_0, t_f]$, $x(t) \in \mathbb{R}^{n_x}$ is the state, $\pi_\theta(x): \mathbb{R}^{n_x} \to \mathbb{R}^{n_u}$ is the state feedback controller and $\theta \in \mathbb{R}^{n_\theta}$ are its parameters.

\section{Methodology}

The most common parametrization in trajectory optimization methods utilizes low-order polynomials with a predefined set of time intervals to approximate the controller as a function of time \citep{rao_SurveyNumerical_2009, teo_AppliedComputational_2021}.  
In contrast, in our work the the control function is posed the output of a parametrized state feedback controller.
The controller parameters are constant through all the integration time (they define statically the nonlinear feedback controller) and are the only optimization variables of the problem.
This transforms the original dynamical optimization problem into a parameter estimation one \citep{teo_AppliedComputational_2021}.
The approach approximates an optimal closed-loop policy, which is a continuous function that is well-defined for states outside the optimal path.
This quality allows to calculate the phase portrait of the system with the embedded control policy, as shown in section \ref{sec:case_studies}.

\footnotetext[3]{
    Evaluation of $H_x$ and $H_\theta$ during integration require $x(t)$ to be available at each evaluation time.
    It is possible to either store an interpolation of the path, recalculate it from checkpoints on-demand or integrate it again but backwards in time from $x(t_f)$ jointly with the adjoint equations.
    Several works study the possible variations on this point and the trade-offs between accuracy, memory consumption and efficiency\citep{serban_cvodes_2005, ma_ComparisonAutomatic_2021, zhuang_AdaptiveCheckpoint_2020, daulbaev_InterpolationTechnique_2020}.  
    We used the interpolation scheme described in \citep{rackauckas2020universal} due to better numerical stability.
}
\footnotetext[4]{
    The cost gradient can also be computed through quadrature routines on $H_\theta$ instead once both $x(t)$ and $\lambda(t)$ have been integrated and their interpolation is available for $t \in [t_0, t_f]$.
}

In this work, an adaptive step size ODE solver is the procedure that selects the adequate time discretization during the differential equation integration.
The trade-off between computational speed and accuracy is also exclusively set by the adaptive integration method and its error tolerance \citep{DifferentialEquations.jl-2017}.
This allows the integrator to account for the influence of the policy to keep a constant integration error throughout the time interval.

Once the gradient of the objective function is computed, any first-order gradient based optimization procedure can be used.
Through this work we used the ADAMW optimizer \citep{kingma_AdamMethod_2015, loshchilov_DecoupledWeight_2019} for the preconditioning procedure.
The IPOPT optimizer \citep{wachter_ImplementationInteriorpoint_2006} with the limited-memory quasi-Newton approximation was used throughout the Fiacco-McCormick iterations.

\subsection{Adjoint sensitivity analysis} \label{sec:adjoint_sensitivity}

The Euler-Lagrange equations for the Bolza problem are the backbone of the calculation of the gradients of a functional objective with respect to the policy parameters.
Here we display the relevant equations for clarity and describe the precise algorithm for their numerical evaluation using automatic differentiation.

First we consider the augmented objective function with smooth Lagrange multipliers $\lambda(t) \in \mathbb{R}^{n_x}$ (also called costate/adjoint vector) to enforce the dynamical equations as constraints,
\begin{equation*}
    L(x, \theta, \lambda) = J(x, \theta) - \int_{t_0}^{t_f}\lambda^T(\dot x - f(x,\theta)) \, dt.
\end{equation*}
It is convenient to define the Hamiltonian function
\begin{equation*}
    H(x, \theta, \lambda) = \lambda^T f(x, \theta) + \ell(x,\theta).
\end{equation*}
Denoting partial derivatives with subscripts, the dynamic constraints are then equivalent to
\begin{equation} \label{eq:dynamic_equations}
    \dot x(t) = H_\lambda(x, \theta, \lambda).
\end{equation}
The following form of multipliers $\lambda$ assure the stationary path required for first-order optimality ($\mathrm{d} L = 0$):
\begin{subequations} \label{eq:adjoint_ode}
    \begin{align}
        \dot \lambda(t) = - H_x & = -\lambda^T f_x - \ell_x, \\
        \lambda(t_f)            & = \phi_x(x(t_f)).
    \end{align}
\end{subequations}
The sensitivity of the objective function $J$ with respect to the parameters $\theta$ are
\begin{equation} \label{eq:cost_gradient}
    J_\theta = L_\theta = \int_{t_0}^{t_f} \lambda^T f_\theta + \ell_\theta \, dt
    = \int_{t_0}^{t_f} H_\theta \, dt.
\end{equation}
It is important to note that both the adjoint \eqref{eq:adjoint_ode} and the cost sensitivity equations \eqref{eq:cost_gradient} require the availability of the forward path solution $x(t)$ (see \ref{fnote:terminology}).

As in control vector iteration \citep{luus_ControlVector_2009}, the last optimality condition is relaxed and left to be approximated iteratively by the optimization procedure over the controller parameters $\theta$:
\begin{equation} \label{eq:optimal_theta}
    H_\theta = \lambda^T f_\theta + \ell_\theta \rightarrow 0.
\end{equation}
The terms of the form $\lambda^T f_\square$ are vector-Jacobian products that can be efficiently calculated with reverse AD \citep{chen_NeuralOrdinary_2018}.
The precise numerical implementation is shown in algorithm \ref{alg:adjoint_sensitivity}.

The parametrization variables of the controller $\theta$ may be thought of as unconstrained time-constant variables; their optimization is a parameter estimation problem with dynamic constraints \citep*{teo_AppliedComputational_2021}.
Since the optimization variables are the unconstrained parameters of the policy, it is not necessary to use the more general Pontryagin's Maximum Principle, which explicitly deal with control variable constraints.
Further problem variants including minimum-time problems could also be considered without great modifications \citep{bryson_applied_1975}.  

\begin{algorithm} \caption{Adjoint Sensitivity pseudocode}\label{alg:adjoint_sensitivity}
    \begin{algorithmic}[0]
        \State \textbf{Inputs:} Set of parameters $\theta$, neural controller $\pi_\theta$, dynamical system $\dot{x}(t) = f(x,\theta)$, initial condition $x_0$, fixed time interval $[t_0, t_f]$, running cost function $\ell\left(x(t), \pi_\theta(x)\right)$, final cost function $\phi(x)$.
        \Procedure{Forward integration}{}
        \State Calculate $x(t)$ by integrating forward the dynamical system.
        \State Store $x(t_f)$ and an interpolation of the path $x(t)$, or setup checkpoints to recalculate it on demand.\footnotemark
        \EndProcedure
        \State Set the adjoint ODE: $\dot{\lambda} = -H_x$ with boundary condition $\lambda(t_f) = \phi_x(x(t_f))$. \Comment{Equation \ref{eq:adjoint_ode}.}
        \State Set the quadrature ODE: $\dot{J_\theta} = H_\theta$ with boundary condition $J_\theta(t_f) = 0$. \Comment{Equation \ref{eq:cost_gradient}.}
        \Procedure{Backward integration}{}
        \State Define the extended ODE system $\dot{F} \coloneqq \left[-H_x, -H_\theta \right]$.
        \State Integrate $F$ backwards from $t_f$ to $t_0$.
        \EndProcedure
        \State \Output Last block of the extended system at the origin: $F(t_0)$.\footnotemark \Comment{$J_\theta = \int_{t_f}^{t_0} -H_\theta \, dt$}
    \end{algorithmic}
\end{algorithm}

\subsection{Neural Policy}

Given the initial condition $x_0$ and a set of parameters $\theta$, the evolution of the states $x(t),\ \forall t \in [t_0, t_f]$ is obtained by the numerical integration of the dynamical system following the parametrized control $u(t) = \pi_\theta(x(t))$.
In practice, at each point where the system requires to be evaluated by the ODE  integrator, the controller will also be evaluated at the current state and its output sets the control component of the system dynamics (see section \ref{sec:case_studies} for specific examples).

We use a dense neural network with parameters $\theta$ as a state feedback control policy.
The controller input is the state of the system at the current time $x(t)$ and its outputs are the constrained controls  $u(t) \in U \subset \mathbb{R}^{n_u}$,
\begin{equation} \label{eq:control_parametrization}
    \pi_\theta(t) = \pi(x(t), \theta) \in U.
\end{equation}
The output layer of the network always has the same number of nodes as the control dimension.

\subsection{Control constraints}

The controller outputs are box-constrained to $U$ directly from the nonlinearity of the last activation function of the policy.
We use a custom scaled sigmoid function $\sigma(\cdot): \mathbb{R} \to [0,1]$ for each control component $u_i$.
Then, the last layer of the policy has the functional form
\begin{equation*}
    u_i = u_i^{(lb)} + (u_i^{(ub)} - u_i^{(lb)}) \sigma_i(\cdot),
\end{equation*}
where $u_i^{(lb)}$ and $u_i^{(ub)}$ are the lower and upper bounds of the control component, respectively, and $i \in [1,\ldots, n_u]$.
By construction, any parameter set in the policy yields a feasible control sequence satisfying the control constraints.

\subsection{State constraints}
State constraints are enforced via penalty functions added to the original functional cost function as a running cost,
\begin{equation*}
    \ell(x, \pi_\theta) \rightarrow \ell(x, \pi_\theta) + \alpha P(x, \delta),
\end{equation*}
where the weight parameter $\alpha$ and the relaxation parameter $\delta$ are both positive.
This approach purposely avoids the treatment of state constraints through additional Lagrange multipliers.
This would be problematic due to alternating optimality conditions along subarcs with active constraints, unknown in number and locations \textit{a priori} \citep{bryson_applied_1975, biegler_NonlinearProgramming_2010}.  

\subsubsection{Relaxed logarithmic penalties}

\begin{figure}[t]
    \centering
    \includegraphics[width=0.55\columnwidth]{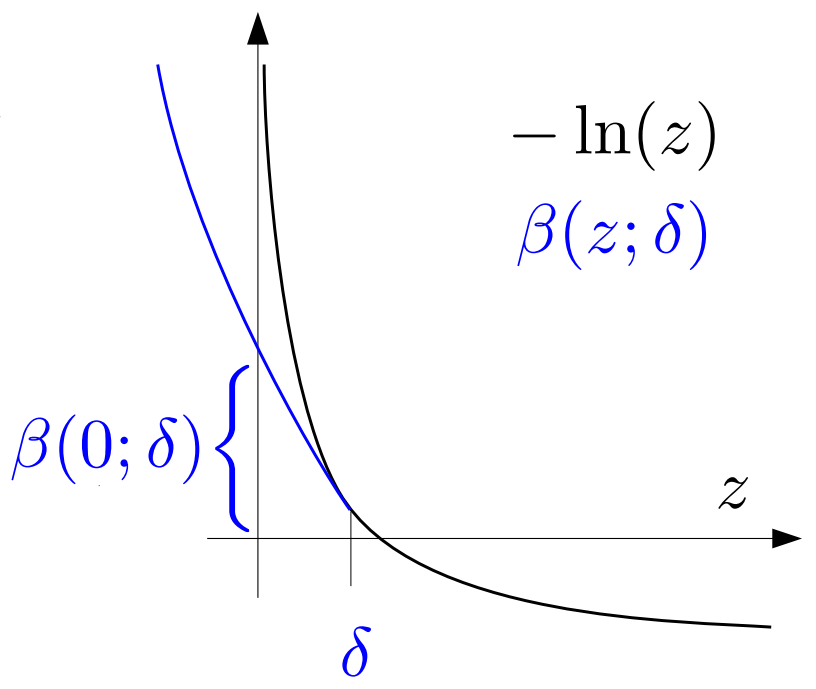}
    \caption{
        Relaxed logarithmic penalty function \citep{feller_ContinuoustimeLinear_2014}.
        Both coefficient and relaxation parameter can be adjusted to approximate interior point penalties iteratively.
    }
    \label{fig:relaxed_log}
\end{figure}

It is not trivial to understand the black-box relationship between the policy parameters  and the specific control profiles the corresponding policy generates.
This is specially problematic because unoptimized parameters may result in a state evolution that violates state constraints.
Consequentially, a state constraint penalty must be well defined for states outside the feasible region to correct the initial profiles.
On the other hand, penalties defined exclusively outside the feasible region lack correction signals in the vicinity of the boundary and inside the allowed control region.
To address this, we leverage a logarithmic relaxation that extends the penalty definition outside the feasible region while preserving an interior feedback as logarithmic barriers used in interior point methods.

We use a variant of the relaxed logarithmic barrier function introduced in \citep{feller_ContinuoustimeLinear_2014} (see figure \ref{fig:relaxed_log}):  
\begin{equation} \label{eq:relaxed_penalty}
    P(z)=\left\{
    \begin{aligned}
        -\ln (z)          &  & z>\delta,       \\
        \beta(z ; \delta) &  & z \leq \delta,
    \end{aligned}
    \right.
\end{equation}
with the truncated exponential relaxing function
\begin{equation*}
    \beta(z; \delta)=\exp (1-z/\delta)-1-\ln (\delta).
\end{equation*}
The scalar $\delta \ge 0$ is called the relaxation parameter.
and controls the distance from the constraint boundary where the relaxation of the barrier function starts.
Note that in equation \ref{eq:relaxed_penalty} the transition between both functions at $\delta$ is smooth.

The progression of barrier parameters follows \textit{Fiacco-McCormick iterations} \citep{fiacco_NonlinearProgramming_1990} with a barrier update after several optimization iterations.
The penalties are increased by tightening both the coefficient ($alpha$) and the relaxation parameter ($delta$) sequentially between optimization rounds where they are left fixed.  
Once the initial parameters ($\theta$) are settled, the initial penalty parameters are tuned to start the first optimization round with a reasonable ratio between the original objective function and the state penalty.
After each optimization procedure, the relaxation parameter $\delta$ is decreased by a constant rate until this ratio is attained again, to then proceed with the next optimization round.
If at any round this reduction becomes problematic for the optimization due to the steep convergence to a logarithmic penalty, the coefficient $\alpha$ is increased instead in that occasion.
The iterations are repeated until the loss function stops decreasing or a predefined number of barrier updates is reached.

\subsection{Policy parametrization}

For setting the initial values of the parameters of the controller we use the Glorot uniform initialization \citep{glorot_UnderstandingDifficulty_2010}: biases are initialized to zero and the weights $W^{(i)}$ of layer $i$ follow a uniform distribution
\begin{equation*}
    W^{(i)} \sim \mathcal{U}\left[-b, b\right], \quad b = G\sqrt{\frac{6}{n_{in} + n_{out}}},
\end{equation*}
where $n_{in}$ and $n_{out}$ are the number of input and output units in that layer, respectively.
The gain $G$ is set to $1$ for the output layer with sigmoid activations and $5/3$ for hidden layers with tanh activations, as is commonly done in practice for this nonlinear activations.
Crucially for our application, even though the Glorot initialization avoids backpropagation pathologies, the control profile that the policy produce is not necessarily a good initial guess for the optimization procedure.

\subsection{Preconditioning}
In this section we describe several strategies that proved useful to improve the convergence of the optimization to a feasible policy.
Similarly to indirect methods in optimal control \citep{rao_SurveyNumerical_2009}, the convergence improves when the initial parameters produces a good initial control profile.
Thus, it is beneficial to optimize towards a reasonable control profile as a preconditioning step.
Afterwards, the original problem can be optimized starting from the policy parameters that reproduce the preconditioning control profile.

The optimization towards a predefined profile is a special case of the problem statement (equation \ref{eq:problem_statement}), with a running cost that penalizes the discrepancy from the profile along the time evolution of the system,
\begin{equation*}
    J_{pre} = \int_{t_0}^{t_f}{\left[\pi_\theta(x(t)) - u_{ref}(t)\right]^2}\, dt.
\end{equation*}

Since it is challenging to link the network parameter values to the control profile that the policy generates with them, initial stiffness of the policy can also be problematic.
This is alleviated by applying the preconditioning procedure iteratively on gradually increasing time windows, to eventually cover the entire time span.

Moreover, a multi-start method is advisable given the high nonlinearity of the optimization problem.
Generating multiple feasible control profiles instead of multiple aleatoric network parameters proved more effective to guide the initial control search.

Furthermore, solving a coarse version of the same optimal control problem (without policy) through classic methods like collocation \citep{rao_SurveyNumerical_2009} provides the best initial guess as a function of time $u(t)$, since it is closer to the optimal solution.
A related approach that also relies on a collocation solution to guide the NeuralODE training was suggested in \citep{roesch_CollocationBased_2021a}, yet applied to identify mechanistic models in a data-driven manner instead of control problems.
The collocation routines used Gauss-Lobatto quadrature with low-order polynomials over a coarse time grid covering the same time span of the original problem \citep{pulsipher_UnifyingModeling_2022}.

\section{Case studies} \label{sec:case_studies}

We applied the discussed techniques to a set of optimal control problems with nonlinear dynamical systems \footnote{These case studies are particular instances of the general problem statement in equation \ref{eq:problem_statement}.}.
Since the state feedback controller is well defined in continuous time and space, it is possible to trace the whole phase space of the dynamical system with the embedded controller.
We show the convergence of the policy towards feasible solutions satisfying constraints through the iterative tightening of the relaxed logarithmic penalties.

\subsection{Van der Pol Oscillator}

\begin{figure*}[t]
    \centering
    \subfloat{\includegraphics[width=\columnwidth]{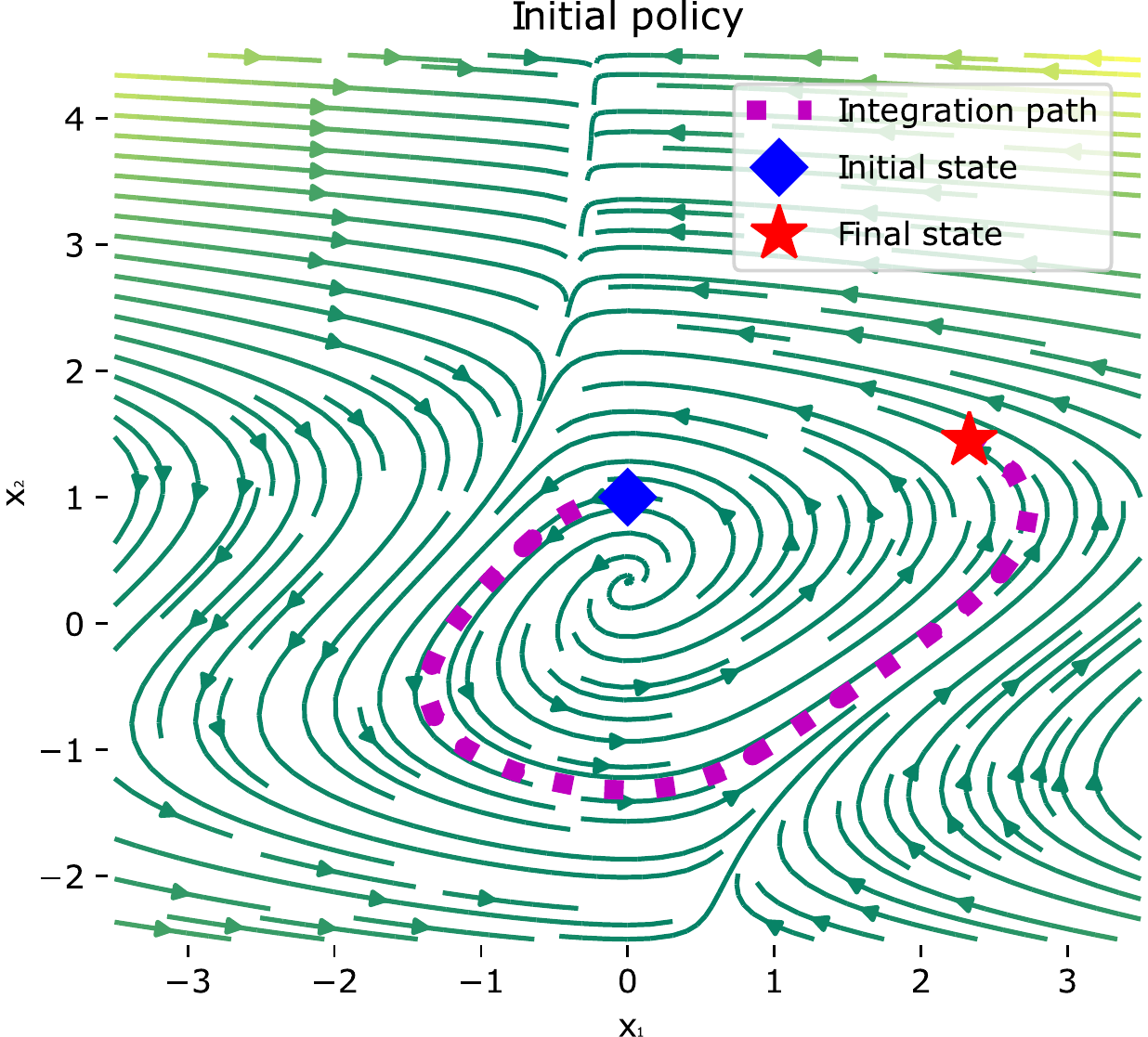}
        \label{fig:vdp_phase_initial}}
    \hfil
    \subfloat{\includegraphics[width=\columnwidth]{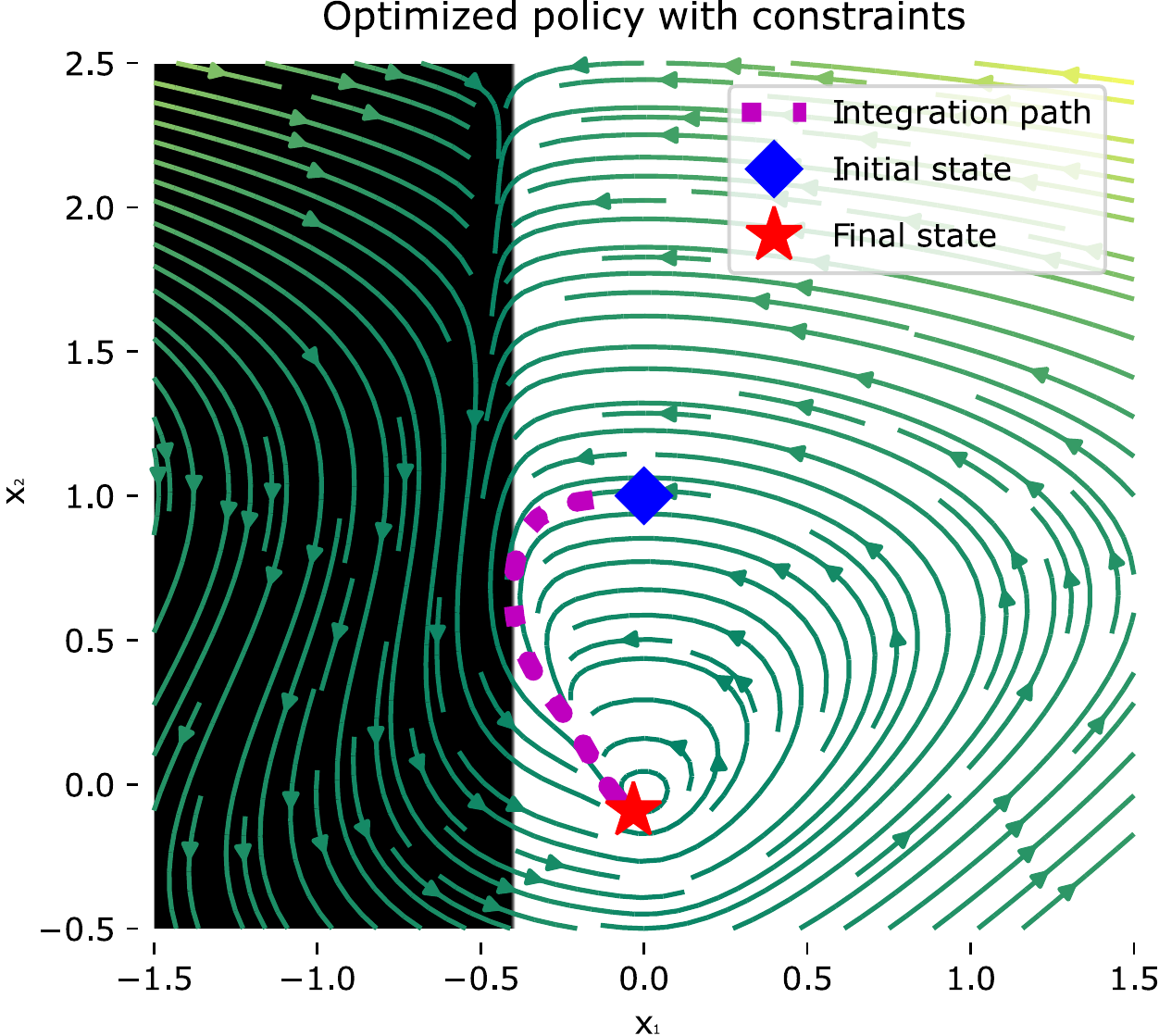}
        \label{fig:vdp_phase_optimized_constrained}}
    \caption{
        Phase portraits of the Van der Pol system with the embedded neural controller.
        The portrait on the left corresponds to an unoptimized policy while the portrait on the right corresponds to the optimized policy.
        The transition between portraits shows how the system dynamics change after the controller has been optimized and complies with the state constraints (shown as the black region).
    }
    \label{fig:vdp_phase_portrait}
\end{figure*}

We analysed the constrained van der Pol oscillator problem from \citep{vassiliadis_SolutionClass_1994}:
\begin{equation} \label{eq:van_der_pol}
    \begin{aligned}
        \min_{\theta} \quad & J = \int_0^5 x_1^2 + x_2^2 + u ^2 \,dt, \\
        \textrm{s.t.} \quad & \dot x_1(t) = x_1(1-x_2^2) - x_2 + u,   \\
                            & \dot x_2 (t) = x_1,                     \\
                            & x_1(t) + 0.4 \geq 0,                    \\
                            & -0.3 \leq u(t) \leq 1,                  \\
    \end{aligned}
\end{equation}
with initial condition $x(t_0) = (0, 1)$.
The state feedback controller was a dense neural network with 2 hidden layers of 16 nodes
\begin{equation*}
    u(t) = \pi_\theta(x_1(t), x_2(t)).
\end{equation*}
Intermediate layers used a hyperbolic tangent while the output layer used a sigmoid scaled to the allowed interval.

The problem was first optimized disregarding the running constraints (figure \ref{fig:vdp_phase_initial}), obtaining a final objective $J$ of 2.87.
This optimal solution was subsequently used as the starting parameters for the constrained version of the problem, where the running penalty was added to enforce state constraints.
The final constrained solution is displayed in figure \ref{fig:vdp_phase_optimized_constrained}, where the enforced constraint on $x_1(t)$ is partially active along the optimal path.
The optimum value was 2.953 in the constrained case; both optimal values match closely the results reported in \citep{vassiliadis_SolutionClass_1994}.

Since the policy is embedded within the dynamical system as a state feedback controller, the policy is well defined for states outside the optimal trajectory and modifies the whole phase space picture.

\subsection{Bioreactor}

As a more challenging case study of chemical engineering relevance \citep{vassiliadis2021optimization}, we optimize the bioreactor model analysed in \citep{bradford_StochasticDatadriven_2020}.
The reactor has 3 state variables and 2 control functions, which represent the biomass concentration ($C_X$), the nitrate concentration ($C_N$), the bioproduct concentration ($C_{q_{c}}$), the light intensity ($I$) and the nitrate inflow rate ($F_N$) of the system, respectively:  

\begin{equation*}
    \begin{aligned}
        \frac{d C_{X}}{d t}     & =\frac{u_{m} C_{X} I}{I+k_{s}+\frac{I^{2}}{k_{i}}} \frac{C_{N}}{C_{N}+K_{N}}-u_{d} C_{X},             \\
        \frac{d C_{N}}{d t}     & =-Y_{\frac{N}{X}}  \frac{u_{m} C_{X} I}{I+k_{s}+\frac{I^{2}}{k_{i}}} \frac{C_{N}}{C_{N}+K_{N}}+F_{N}, \\
        \frac{d C_{q_{c}}}{d t} & = \frac{k_{m} C_{X} I}{I+k_{s q}+\frac{I^{2}}{k_{i q}}} - \frac{k_{d} C_{q_{c}}}{C_{N}+K_{N p}},
    \end{aligned}
\end{equation*}
with the initial condition
\begin{equation*}
    [C_X(0), C_N(0), C_{q_c}(0)] = [1, 150, 0].
\end{equation*}
The rest of the parameters are arranged in table \ref{tab:bioreactor} in the appendix.

The main objective is to maximize the bioproduct $C_{q_c}$ at the end of the batch time $T$.
Additionally, there are 2 running state constraints and 1 final state constraint conditions that must be satisfied, as well as 2 control constraints:

\begin{equation*}
    \begin{aligned}
        C_{N} \le 800,                 & \qquad \forall t \in (0, t_f), \\
        0.011 C_{X} - C_{q_c} \le 0.3, & \qquad \forall t \in (0, t_f), \\
        C_{N} \le 150,                 & \qquad t=t_f,                  \\
        120 \le I \le 400,             & \qquad \forall t \in (0, t_f), \\
        0 \le F_N \le 40,              & \qquad \forall t \in (0, t_f).
    \end{aligned}
\end{equation*}

As before, the controller has 2 hidden layers of 16 nodes, with 3 input and 2 output nodes
\begin{equation*}
    I(t), F_N(t) = \pi_\theta(C_X(t), C_N(t), C_{q_c}(t)).
\end{equation*}

\begin{figure*}[t]
    \centering
    \subfloat[Running and final constraints.]{\includegraphics[width=\columnwidth]{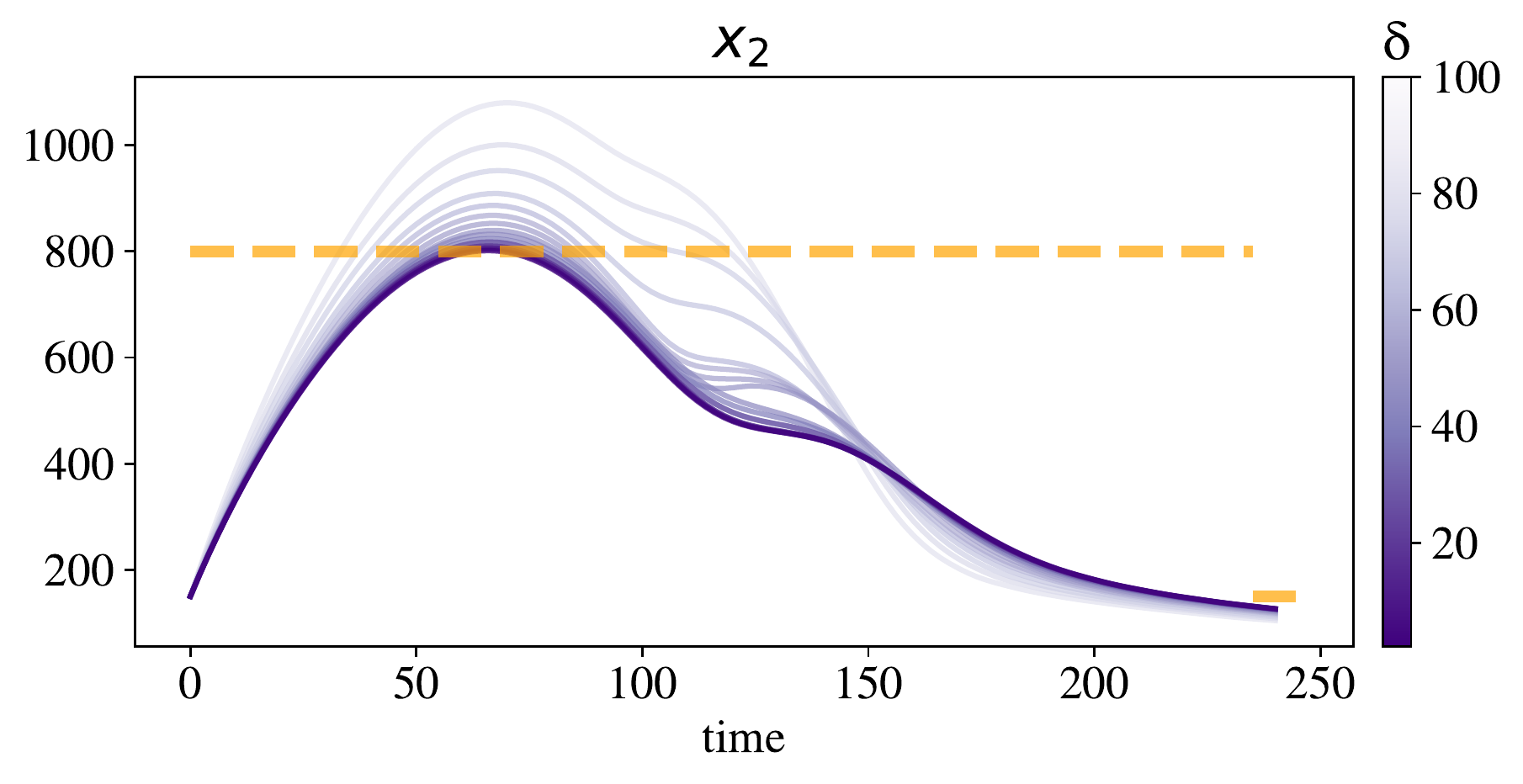}%
        \label{fig:bioreactor_constraint_x2}}
    \hfil
    \subfloat[Running compound constraint.]{\includegraphics[width=\columnwidth]{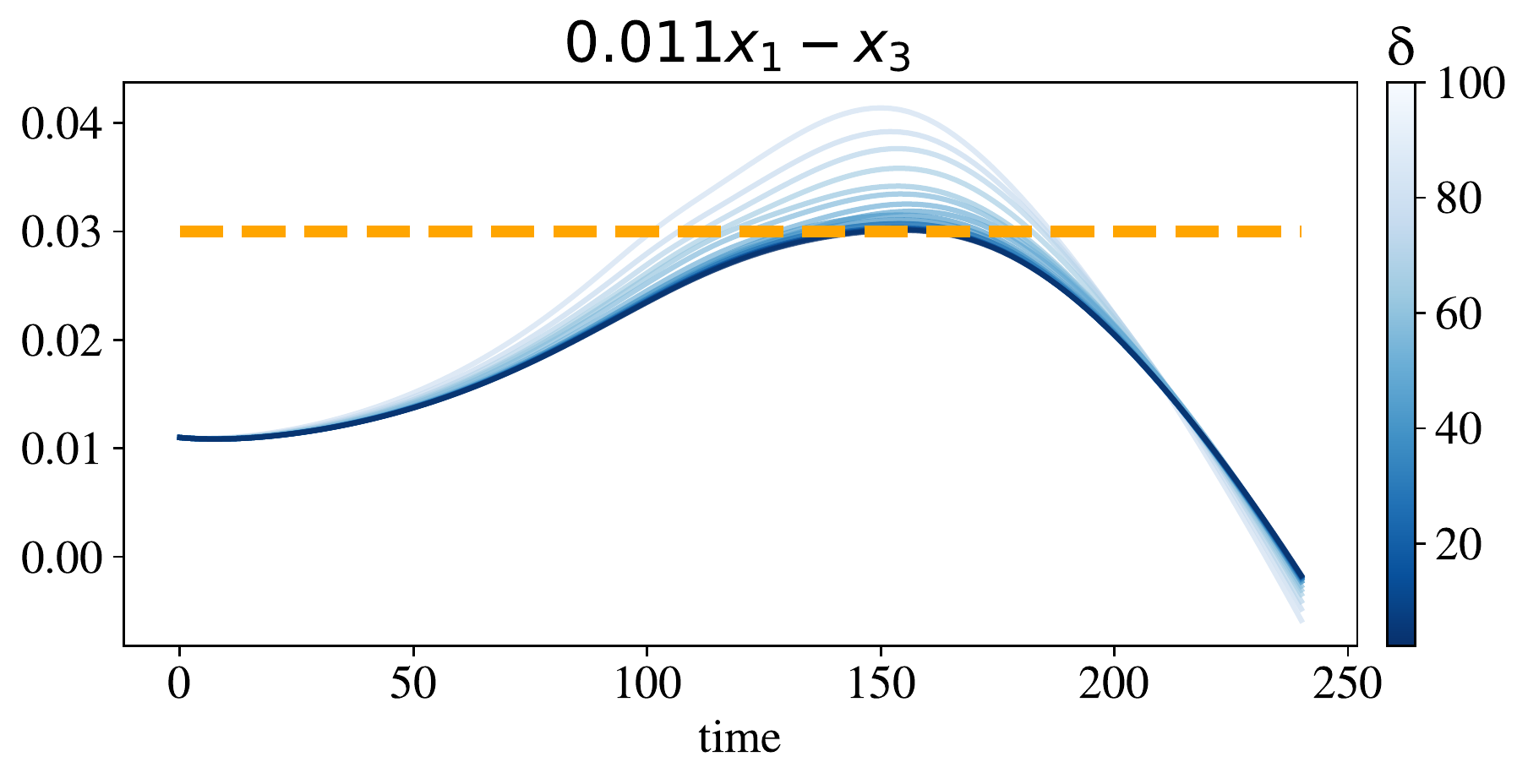}%
        \label{fig:bioreactor_constraint_x1_x3}}
    \caption{
        Constraint enforcement through relaxed logarithmic penalty tightening.
        As the $\delta$ parameter is reduced, penalties tend to the classic logarithmic penalties from interior point methods.
        Higher values of $delta$ are useful to deal with infeasible paths generated by unoptimized policies in the initial iterations.
    }
    \label{fig:bioreactor_constraints}
\end{figure*}

\begin{table}[hb]
    \renewcommand{\arraystretch}{1.2}
    \caption{Bioreactor parameters \citep{bradford_StochasticDatadriven_2020}.}
    \label{tab:bioreactor}
    \centering
    \begin{tabular}{cc}
        \hline
        Parameter             & Value   \\ \hline
        $u_m$                 & 0.0572  \\
        $u_{d}$               & 0       \\
        $K_{N}$               & 393.1   \\
        $Y_{\nicefrac{N}{X}}$ & 504.5   \\
        $k_{m}$               & 0.00016 \\
        $k_{d}$               & 0.281   \\
        $k_{s}$               & 178.9   \\
        $k_{i}$               & 447.1   \\
        $k_{sq}$              & 23.51   \\
        $k_{iq}$              & 800     \\
        $K_{N_p}$             & 16.89   \\ \hline
    \end{tabular}
\end{table}

\section{Conclusion}
In this work, we explored how neural feedback policies within a white-box dynamical model may be used to solve nonlinear optimal control problems in continuous time.
We also showcased practical approaches to satisfy both state and control constraints with the optimized policy, which are crucial for real world applications and an active area of research.
The approach takes motifs from reinforcement learning by leveraging a neural network as a nonlinear controller, within an environment defined by a deterministic dynamical system.
We relied on adjoint sensitivity analysis to handle the large number of parameters required by the neural policy.
This approach takes advantage of the specific ODE model form for a precise gradient calculation that avoids sampling completely.
In contrast to open-loop methods in optimal control, the state feedback policy can be evaluated for any state outside the optimal trajectory.
Our policies provide the whole phase space of the controlled nonlinear system, approximating the intractable closed-loop solution of the nonlinear optimal control problem.

\begin{ack}
    The authors would like to thank Christopher Rackauckas and the Julia language community for helpful discussions on the numerical implementation of this work.
    We are grateful to Miguel de Carvalho and Fabian Thiemann for their helpful reviews and comments.
\end{ack}

\bibliography{ifacconf}             

\end{document}